\documentclass[leqno]{amsart}

\usepackage{latexsym}
\def\qed{\hfill $\Box$} 
\def\tr#1{\mathord{\mathopen{{\vphantom{#1}}^t}\!#1}}
\usepackage{amsthm}
\usepackage{graphicx,color}
\usepackage{amsmath}
\usepackage{amsfonts}
\usepackage{exscale}
\usepackage{enumerate}

\newtheorem{theo}{Theorem}[section]
\newtheorem{defi}{Definition}[section]
\newtheorem{lemm}{Lemma}[section]
\newtheorem{coro}{Corollary}[section]

\newtheorem{pro}{Proposition}[section]
\makeatletter
\@addtoreset{equation}{subsection}

\newfont{\bg}{cmr9 scaled\magstep4}
\newcommand{\bigzerol}{\smash{\lower1.0ex\hbox{\bg 0}}}
\newcommand{\bigzerou}{
\smash{\hbox{\bg 0}}}
\newcommand{\R}{\mathbb{R}}

\usepackage{amssymb}

\title{Distance-squared mappings}

\author{Shunsuke Ichiki and   
Takashi Nishimura}
\address{Graduate School of Environment and Information Sciences,
Yokohama National University,
Yokohama 240-8501, Japan.}
\email{ichiki-shunsuke-jb@ynu.ac.jp}
\address{Research Group of Mathematical Sciences,
Research Institute of Environment and Information Sciences,
Yokohama National University,
Yokohama 240-8501, Japan.}
\email{nishimura-takashi-yx@ynu.jp}


\begin{document}
\begin{abstract}
A distance-squared function is one of the most significant functions 
in the application of singularity theory to differential geometry.
In this paper, we define naturally extended mappings of distance-squared functions, wherein each component is a distance-squared function. 
We investigate the properties of these mappings from the viewpoint of differential topology. 
\end{abstract}
\date{}
\maketitle
\noindent
\section{Introduction}\label{IN}

Let $\mathbb{R}^n$ be the $n$-dimensional Euclidean space; thus, 
a point $x\in \mathbb{R}^n$ is $n$-tuple $x=(x_1,\ldots,x_n)$ of real numbers.
Let $d : \mathbb{R}^n\times \mathbb{R}^n\rightarrow \mathbb{R}$ be the 
$n$-dimensional Euclidean distance,
\[d(x,y)=\sqrt{\sum_{i=1}^n(x_i-y_i)^2},\]
where $x=(x_1,\ldots,x_n)$ and $y=(y_1,\ldots,y_n)$.
Let $p$ be a given point in $\mathbb{R}^n$. 
The mapping 
$d_p : \mathbb{R}^n\rightarrow \mathbb{R}$, defined by $d_p(x)=d(p,x)$, is called a {\it distance function}.
\begin{defi}
{\rm Let $p_1,\ldots,p_\ell$ $(\ell\geq 1)$ be $\ell$ given points in $\mathbb{R}^n$. 
The mapping 
$d_{(p_1,\ldots,p_\ell)} : \mathbb{R}^n \rightarrow \mathbb{R}^\ell$ defined by 
\[
d_{(p_1,\ldots,p_\ell)}(x)=(d(p_1,x),\ldots,d(p_\ell,x)), 
\] 
is called a {\it distance mapping}.}
\end{defi}
A distance mapping is one in which each component is a distance function.

Let $S^n$ be the $n$-dimensional unit sphere in $\mathbb{R}^{n+1}$.   
\begin{pro} \label{d}
Let $i : S^1\rightarrow i(S^1)\subset \mathbb{R}^2$ be a homeomorphism.
Then, there exist two points $p_1$,$p_2$ $\in i(S^1)$ such that $d_{(p_1,p_2)}\circ i $ is a homeomorphism to 
its image $\left(d_{(p_1,p_2)}\circ i\right)(S^1)$.
\end{pro}
\begin{figure}
\begin{center}
\includegraphics[width=6cm,clip]{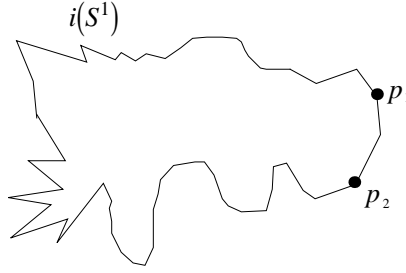}
\caption{Figure of a mapping $i: S^1\to \mathbb{R}^2$ for which Proposition \ref{d} can be applied.}
\end{center}
\end{figure}\label{figure 1}
Proposition \ref{d} is proved in Section \ref{dP}.
Although Proposition \ref{d} is applicable even if a mapping $i$ is not differentiable 
(see Figure \ref{figure 1}), it appears difficult to derive higher-dimensional extensions of the proposition.   

On the other hand, the differentiable version of higher-dimensional extensions can be obtained as follows. 
All manifolds and mappings in this paper belong to class $C^\infty $ except in Proposition \ref{d} 
and the homeomorphism $\sqrt{\quad}: \mathbb{R}^\ell_+\to \mathbb{R}^\ell_+$ defined below. 
\begin{defi}
{\rm Let $p_1,\ldots,p_\ell$ $(\ell\geq 1)$ be $\ell$ given points in $\mathbb{R}^n$. 
The mapping $D_{(p_1,\ldots,p_\ell)} : \mathbb{R}^n \rightarrow \mathbb{R}^\ell$ 
defined by 
\[
D_{(p_1,\ldots,p_\ell)}(x)=(d^2(p_1,x),\ldots,d^2(p_\ell,x)),
\]
is called a {\it distance-squared mapping}}.   
\end{defi}
Note that $D_{(p_1,\ldots,p_\ell)}$ always has a singular point if $\ell \le n$.       
Nevertheless, Theorems \ref{D1} and \ref{D2} below hold.      
\begin{theo} \label{D1}
Let $M$ be an $m$-dimensional closed manifold $(m\geq 1)$, and let $i : M\to \mathbb{R}^\ell$ 
$(m+1\leq \ell)$ be an embedding. Then, there exist $p_1,\ldots,p_{m+1}\in i(M)$, $p_{m+2},\ldots,p_\ell\in \mathbb{R}^\ell $such that $D_{(p_1,\ldots,p_\ell)}\circ i : M\rightarrow \mathbb{R}^\ell$ is an embedding.   
\end{theo}
\begin{coro}
Let $M$ be an $m$-dimensional closed manifold $(m\geq 1)$, and let $i : M\rightarrow \mathbb{R}^{m+1}$ be an embedding. Then, there exist $p_1,\ldots,p_{m+1}\in i(M)$ 
such that $D_{(p_1,\ldots,p_{m+1})}\circ i : M\rightarrow \mathbb{R}^\ell$ is an embedding.  
\end{coro}   
Define $\mathbb{R}_+^\ell=\{(x_1, \ldots, x_\ell)\; |\; x_i\ge 0\; (1\le i\le \ell)\}$ and let 
$\sqrt{\quad} : \mathbb{R}_+^\ell\to \mathbb{R}_+^\ell$ be the homeomorphism defined by 
$\sqrt{\quad}(x_1, \ldots, x_\ell)=(\sqrt{x_1}, \ldots, \sqrt{x_\ell})$.   It is clear that 
$d_{(p_1,\ldots,p_\ell)}=\sqrt{\quad}\circ D_{(p_1,\ldots,p_\ell)}$.    Thus, the following corollary holds: 
\begin{coro}
Let $M$ be an $m$-dimensional closed manifold $(m\geq 1)$, and let $i : M\to \mathbb{R}^\ell$ 
$(m+1\leq \ell)$ be an embedding. 
Then, there exist $p_1,\ldots,p_{m+1}\in i(M)$, $p_{m+2},\ldots,p_\ell\in \mathbb{R}^\ell $such that $d_{(p_1,\ldots,p_\ell)}\circ i : M\rightarrow 
\left(d_{(p_1,\ldots,p_\ell)}\circ i\right)(M)\subset \mathbb{R}^\ell$ is a homeomorphism.     
\end{coro} \label{corollary 1.2}

Let $M$ and $N$ be manifolds. An immersion $f : M\rightarrow N$ 
is said to be 
{\it with normal crossings at a point $y\in N$} if $f^{-1}(y)$ is a finite set $\{x_1,x_2,\ldots,x_k\}$ and for any subsets $\{\lambda _1,\lambda _2,\ldots,\lambda _s\}\subset 
\{1,2,\ldots,k\}$ $(s\leq k)$, 
\[{\rm codim}\,\Bigl(\bigcap _{j=1}^sdf_{x_{\lambda _{j}}}(T_{x_{\lambda _{j}}}M)\Bigr)=\sum_{j=1}^s{\rm codim}\,(df_{x_{\lambda _{j}}}(T_{x_{\lambda _{j}}}M)),\]
where ${\rm codim}\,H={\rm dim}\,T_yN-{\rm dim}\,H$ for any given linear subspace $H$ of $T_yN$. 
An immersion $f : M\rightarrow N$ is said to be {\it with normal crossings} if $f$ is a mapping with normal crossings 
at any point $y\in N$.
\begin{theo} \label{D2}
Let $M$ be an $m$-dimensional closed manifold $(m\geq 1)$, and let $i : M\rightarrow \mathbb{R}^\ell$ $(m+1\leq \ell)$ be an immersion with normal crossings. Then, there exist $p_1,\ldots,p_{m+1}\in i(M)$, $p_{m+2},\ldots,p_\ell\in \mathbb{R}^\ell$
such that $D_{(p_1,\ldots,p_\ell)}\circ i : M\rightarrow \mathbb{R}^\ell$ is an immersion with normal crossings.
\end{theo} 
\begin{coro}
Let $M$ be an $m$-dimensional closed manifold $(m\geq 1)$, and let 
$i : M\rightarrow \mathbb{R}^{m+1}$ be an immersion with normal crossings. Then, there exist $p_1,\ldots,p_{m+1}\in i(M)$ such that $D_{(p_1,\ldots,p_{m+1})}\circ i : M\rightarrow \mathbb{R}^{m+1}$ is an immersion with normal crossings.
\end{coro}

$L$-points $p_1,\ldots,p_\ell\in  \R^n$ $(1\leq \ell\leq n+1)$ are said to be {\it in general position}  if $\ell=1$ or 
$\overrightarrow{p_1p_2},\ldots,\overrightarrow{p_1p_\ell}$ $(2\leq \ell\leq n+1)$ are linearly independent.  
\par 
A mapping $f : \R^n\rightarrow \R^\ell$ is said to be {\it $\mathcal{A}$-equivalent} to a mapping $g : \R^n\rightarrow \R^\ell$ 
if there exist diffeomorphisms $\varphi : \R^n\rightarrow \R^n$ and $\psi : \R^\ell\rightarrow \R^\ell$ 
such that $\psi \circ f\circ \varphi^{-1} =g$.    
\par 
A mapping $f : \R^n\rightarrow \R^\ell$ $(\ell\leq n)$ is called 
{\it the normal form of definite fold mappings} if 
$f(x_1,\ldots,x_n)=(x_1,\ldots,x_{\ell-1}, x^2_\ell+\cdots +x^2_n)$.
\par 
The properties of distance-squared mappings are of significance in the proofs of Theorems 
\ref{D1} and \ref{D2}.
It turns out that distance-squared mappings are $\mathcal{A}$-equivalent to the normal form of 
definite fold mappings if $\ell$-points $p_1,\ldots,p_\ell \in \R^n$ $(2\leq \ell\leq n)$ are in general position 
(see Proposition \ref{P}).
\medskip

In Section \ref{SP}, several preliminaries for the proofs of 
Theorems \ref{D1}, \ref{D2} and Proposition \ref{d} are given.   
The proofs of Theorems \ref{D1} and \ref{D2} are given in Section \ref{SD}. 
Finally, in Section \ref{dP}, the proof of Proposition \ref{d} is given for the sake of readers' convenience.
\section{Preliminaries}\label{SP}
\subsection{Preliminaries for the proofs of Theorems \ref{D1} and \ref{D2}}
\begin{pro}\label{P}
\ \\
$(I)$\ Let $\ell$,$n$ be integers such that $2\leq \ell\leq n$, and let $p_1,\ldots,p_\ell\in \mathbb{R}^n$ be in general position. Then, $D_{(p_1,\ldots,p_\ell)} : \R^n\rightarrow \R^\ell$ is $\mathcal{A}$-equivalent 
to the normal form of definite fold mappings.\\
$(II)$\ Let $\ell$,$n$ be integers such that $1\leq n<\ell $, and let $p_1,\ldots,p_{n+1}\in \mathbb{R}^n$ be in general position. Then, $D_{(p_1,\ldots,p_\ell)} : \R^n\rightarrow \R^\ell$ is $\mathcal{A}$-equivalent to the inclusion $(x_1,\ldots,x_n)\mapsto (x_1,\ldots,x_n,0,\ldots,0)$.
\end{pro}
{\it Proof in the case $\ell=n=2$.}\qquad 
Since the proof requires various elementary coordinate transformations, 
we first prove Proposition \ref{P} in the case $\ell=n=2$ for the sake of clearness.    
By definition, the mapping $D_{(p_1, p_2)}$ has the following form:   
\begin{eqnarray*}
{} & { } & D_{(p_1, p_2)}(x_1, x_2) \\ 
{ } & = & \left((x_1-p_{11})^2+(x_2-p_{12})^2, (x_1-p_{21})^2+(x_2-p_{22})^2)\right) \\ 
{ } & = & \left(x_1^2+x_2^2-2\left( p_{11}x_1+p_{12}x_2\right) +p_{11}^2+p_{12}^2\right. ,  \\ 
{ } & { } & \qquad \qquad 
x_1^2+x_2^2-2\left.\left( p_{21}x_1+p_{22}x_2\right) +p_{21}^2+p_{22}^2\right).    
\end{eqnarray*}
\par 
As Step 1, we remove the quadratic terms in one component.      
This can be done by composing an affine diffeomorphism $h_1$ of the target space.  
We may assume that the composition $h_1\circ D_{(p_1, p_2)}$ has the following form:    
 \begin{eqnarray*}
&{}&(h_1\circ D_{(p_1,p_2)})(x_1,x_2)\\
&=&
\biggl(\sum_{j=1}^2 (p_{2j}-p_{1j})(x_j-p_{1j}), \sum_{j=1}^2(x_j-p_{1j})^2\biggl).   
\end{eqnarray*}   
\par 
As Step 2, we reduce the first component to a linear function.          
This can be done by composing  the affine diffeomorphism 
$h_2(x_1, x_2)=(x_1+p_{11}, x_2+p_{12})$ of the source space.    
The composition has the following form:  
\begin{eqnarray*}
&{}&(h_1\circ D_{(p_1,p_2)}\circ h_2 )(x_1,x_2)\\
&=&
\biggl((p_{21}-p_{11})x_1+(p_{22}-p_{12})x_2, x_1^2+x_2^2 \biggl).   
\end{eqnarray*}   
\par 
As Step 3, we reduce the linear function $(p_{21}-p_{11})x_1+(p_{22}-p_{12})x_2$ to $x_1$.    
Since $p_1, p_2$ are in general position, 
this can be done easily by composing  an affine diffeomorphism 
$h_3$ of the target space and a linear diffeomorphism $h_4$ of the source space.    
We may assume that the composition has the following form:  
\[
(h_3\circ h_1\circ D_{(p_1,p_2)}\circ h_2 \circ h_4)(x_1,x_2)
=
\biggl(x_1, \left( x_1+\alpha x_2\right)^2+x_2^2 \biggl).   
\]
\par 
As Step 4, we complete the square of the second component with respect to the variable $x_2$ as follows:   
\begin{eqnarray*}
\left(x_1+\alpha x_2\right)^2+x_2^2 & = & 
x_1^2+2\alpha x_1 x_2+ \alpha^2x_2^2+x_2^2 \\ 
{ } & = & 
\left(1+\alpha^2\right)\left(x_2+\frac{\alpha}{\left(1+\alpha^2\right)}x_1\right)^2+
\left(1-\frac{\alpha^2}{\left(1+\alpha^2\right)}\right)x_1^2.    
\end{eqnarray*}
\par 
As Step 5, eliminating the term 
$\left(1-\frac{\alpha^2}{\left(1+\alpha^2\right)}\right)x_1^2$ is done 
by composing a nonlinear diffeomorphism 
$h_5$ of the target space.  
\par 
As the final step, the quadratic form  
\[
\left(1+\alpha^2\right)\left(x_2+\frac{\alpha}{\left(1+\alpha^2\right)}x_1\right)^2
\] 
is reduced to $x_2^2$ by composing an affine diffeomorphism $h_6$ of the source space.    
\qed 
\par 
\medskip 
{\it Proof of $(I)$ in general case.}\qquad 
Let $H_1 : \mathbb{R}^\ell\rightarrow \mathbb{R}^\ell$ be the diffeomorphism given by 
\begin{eqnarray*}
{}&{}&H_1(X_1,X_2,\ldots, X_\ell)\\
{}&=&\biggl(\frac{1}{2}\Bigl(X_1-X_2+\sum_{j=1}^n (p_{1j}-p_{2j})^2\Bigl),\ldots, 
\frac{1}{2}\Bigl(X_1-X_\ell+\sum_{j=1}^n (p_{1j}-p_{\ell j})^2\Bigl),X_1\biggl).   
\end{eqnarray*}
The composition of $D_{(p_1,\ldots,p_\ell)}$ and $H_1$ is given by 
\begin{eqnarray*}
&{}&(H_1\circ D_{(p_1,\ldots,p_\ell)})(x_1,x_2,\ldots,x_n)\\
&=&
\biggl(\sum_{j=1}^n (p_{2j}-p_{1j})(x_j-p_{1j}),\ldots,
\sum_{j=1}^n (p_{\ell j}-p_{1j})(x_j-p_{1j}),\sum_{j=1}^n(x_j-p_{1j})^2\biggl).   
\end{eqnarray*}

Let $H_2 : \mathbb{R}^n\rightarrow \mathbb{R}^n$ be the diffeomorphism defined by
\[H_2(x_1,x_2,\ldots,x_n)=(x_1+p_{11},x_2+p_{12},\ldots,x_n+p_{1n}).\]
The composition of $H_1\circ D_{(p_1,\ldots,p_\ell)}$ and $H_2$ is given by 
\begin{eqnarray*}
{}&{}&(H_1\circ D_{(p_1,\ldots,p_\ell)}\circ H_2)(x_1,x_2,\ldots,x_n)\\
{}&=&\biggl(\sum_{j=1}^n (p_{2j}-p_{1j})x_j,\ldots,\sum_{j=1}^n (p_{\ell j}-p_{1j})x_j,\,\sum_{j=1}^nx_j^2\biggl)\\
{}&=&
\begin{pmatrix}
x_1&x_2&\cdots &x_n
\end{pmatrix}
\begin{pmatrix}
p_{21}-p_{11}&\cdots &p_{\ell 1}-p_{11}&x_1\\
p_{22}-p_{12}&\cdots &p_{\ell 2}-p_{12}&x_2\\
\vdots & & \vdots &\vdots \\
\vdots &   & \vdots &\vdots \\
p_{2n}-p_{1n}&\cdots &p_{\ell n}-p_{1n}&x_n
\end{pmatrix}.
\end{eqnarray*}
Set 
\begin{eqnarray*}
A=\begin{pmatrix}
p_{21}-p_{11}&\cdots &p_{\ell 1}-p_{11}\\
p_{22}-p_{12}&\cdots &p_{\ell 2}-p_{12}\\
\vdots & & \vdots \\
\vdots &   & \vdots \\
p_{2n}-p_{1n}&\cdots &p_{\ell n}-p_{1n}
\end{pmatrix}.
\end{eqnarray*}
Since $\ell$-points $p_1,\ldots,p_\ell$ are in general position, it is clear that the rank of $A$ is $\ell -1$.
Hence, by composiong linear coordinate transformations if necessary, we may assume 
that there  exists an $(\ell-1)\times (\ell-1)$ regular matrix $B$ such that
\begin{eqnarray*}
AB=\begin{pmatrix}
  1            &     &      \bigzerol \\
                   & \ddots &       \\
   \bigzerou &     &        1     \\
   \alpha _{11}  & \cdots &\alpha _{1,\ell -1}\\
   \vdots & &\vdots \\
   \alpha _{n-\ell +1,1}  & \cdots &\alpha _{n-\ell +1,\ell -1}
\end{pmatrix}.
\end{eqnarray*}

Let $H_3 : \mathbb{R}^\ell\rightarrow \mathbb{R}^\ell$ be the diffeomorphism defined by
\begin{eqnarray*}
{}&{}&H_3(X_1,X_2,\ldots,X_\ell )=
\begin{pmatrix}
X_1&X_2&\cdots &X_\ell
\end{pmatrix}
\left(
\begin{array}{ccc|c}
&&&  0\\ 
&\mbox{\smash{\Large $B$}}&&\vdots\\
&&&0\\ \hline
0&\cdots &0&1
\end{array}
\right).
\end{eqnarray*}
The composition of $H_3$ and $H_1\circ D_{(p_1,\ldots,p_\ell)}\circ H_2$ is given by 
\begin{eqnarray*}
{}&{}&(H_3\circ H_1\circ D_{(p_1,\ldots,p_\ell)}\circ H_2)(x_1,x_2,\ldots,x_n)\\
{}&=&\begin{pmatrix}
x_1&x_2&\cdots &x_n
\end{pmatrix}
\left(
\begin{array}{ccc|c}
&&&  x_1\\ 
&&&x_2\\
&\mbox{\smash{\Large $A$}}&&\vdots \\
&&&\vdots\\
&&&x_n
\end{array}
\right)
\left(
\begin{array}{ccc|c}
&&&  0\\ 
&\mbox{\smash{\Large $B$}}&&\vdots\\
&&&0\\ \hline
0&\cdots &0&1
\end{array}
\right)\\
{}&=&\begin{pmatrix}
x_1&x_2&\cdots &x_n
\end{pmatrix}
\begin{pmatrix}
  1            &     &      \bigzerol & x_1\\
                   & \ddots &     & \vdots  \\
   \bigzerou &     &        1  & x_{\ell -1}   \\
   \alpha _{11}  & \cdots &\alpha _{1,\ell -1} & x_\ell \\
   \vdots & &\vdots & \vdots \\
   \alpha _{n-\ell +1,1}  & \cdots &\alpha _{n-\ell +1,\ell -1} & x_n
\end{pmatrix}\\
{}&=&\biggl(x_1+\sum_{j=1}^{n-\ell +1}\alpha _{j1}x_{\ell +j-1},\ldots,
x_{\ell -1}+\sum_{j=1}^{n-\ell +1}\alpha _{j,\ell -1}x_{\ell +j-1},\sum_{j=1}^nx_j^2\biggl).   
\end{eqnarray*}
Then, we have
\begin{eqnarray*}
\biggl(x_1+\sum_{j=1}^{n-\ell +1}\alpha _{j1}x_{\ell +j-1},\ldots,
x_{\ell -1}+\sum_{j=1}^{n-\ell +1}\alpha _{j,\ell -1}x_{\ell +j-1},x_\ell ,\ldots,x_n\biggl)\\
=\begin{pmatrix}
x_1&x_2&\cdots &x_n
\end{pmatrix}
\begin{pmatrix}
1     &    &      &     &  & &\\
     &  \ddots  &      &     &  &  &\\
 \bigzerou   && & &\text{\huge{0}}&&\\
 &  & 1 & &&&\\
\alpha _{11}      & \cdots  & \alpha _{1,\ell -1}   &1&&\\
\vdots &       &  \vdots    & 0 & \ddots & &\\
\vdots &       &    \vdots  & \vdots &  & \ddots &\\
\alpha _{n-\ell +1,1}      & \cdots & \alpha _{n-\ell +1,\ell -1}   &0&\cdots &0&1\\
\end{pmatrix}.
\end{eqnarray*}
Set 
\begin{eqnarray*}
C=
\begin{pmatrix}
1     &    &      &     &  & &\\
     &  \ddots  &      &     &  &  &\\
 \bigzerou   && & &\text{\huge{0}}&&\\
 &  & 1 & &&&\\
\alpha _{11}      & \cdots  & \alpha _{1,\ell -1}   &1&&\\
\vdots &       &  \vdots    & 0 & \ddots & &\\
\vdots &       &    \vdots  & \vdots &  & \ddots &\\
\alpha _{n-\ell +1,1}      & \cdots & \alpha _{n-\ell +1,\ell -1}   &0&\cdots &0&1\\
\end{pmatrix}.
\end{eqnarray*}
The inverse matrix of $C$ is the following:
\begin{eqnarray*}
C^{-1}=
\begin{pmatrix}
1     &    &      &     &  & &\\
     &  \ddots  &      &     &  &  &\\
 \bigzerou   && & &\text{\huge{0}}&&\\
 &  & 1 & &&&\\
-\alpha _{11}      & \cdots  & -\alpha _{1,\ell -1}   &1&&\\
\vdots &       &  \vdots    & 0 & \ddots & &\\
\vdots &       &    \vdots  & \vdots &  & \ddots &\\
-\alpha _{n-\ell +1,1}      & \cdots & -\alpha _{n-\ell +1,\ell -1}   &0&\cdots &0&1\\
\end{pmatrix}.
\end{eqnarray*}

Let $H_4 : \R^n\rightarrow \R^n$ be the diffeomorphism defined by 
\[
H_4(x)=xC^{-1}.
\]
The composition of $H_3\circ H_1\circ D_{(p_1,\ldots,p_\ell)}\circ H_2$ and $H_4$ is as follows:
\begin{eqnarray*}
{}&{}&(H_3\circ H_1\circ D_{(p_1,\ldots,p_\ell)}\circ H_2\circ H_4)(x_1,x_2,\ldots,x_n)\\
{}&=&\biggl(x_1,x_2,\ldots,x_{\ell -1},\sum_{k=1}^{\ell -1}\Bigl(x_k-\sum_{j=1}^{n-\ell +1}\alpha _{jk}x_{\ell+ j-1}\Bigl)^2+
\sum_{k=\ell}^nx_k^2\biggl).
\end{eqnarray*}
By completing the square of the last component function 
with respect to variables $x_\ell, \ldots, x_n$, we have the following:   
\begin{eqnarray*}
{ } & { } & 
\sum_{k=1}^{\ell -1}\Bigl(x_k-\sum_{j=1}^{n-\ell +1}\alpha _{jk}x_{\ell+ j-1}\Bigl)^2+\sum_{k=\ell}^nx_k^2 \\ 
{ } & = & 
\sum_{j=\ell}^n\left(\sum_{k=1}^{\ell -1}\alpha_{jk}^2+1\right)
\left(x_j-\frac{1}{\left(\sum_{k=1}^{\ell -1}\alpha_{jk}^2+1\right)}\sum_{k=1}^{\ell-1}\alpha_{j-\ell+1, k}x_k\right)^2 \\ 
{} & {} & \qquad \qquad 
+\sum_{k=1}^{\ell-1}\left(\frac{\sum_{j=\ell}^n\alpha_{j-\ell+1, k}^2}{\left(\sum_{k=1}^{\ell-1}\alpha_{jk}^2+1\right)}
+1\right)x_k^2
\end{eqnarray*}
%
%
\par 
Let $H_5 : \R^\ell\rightarrow \R^\ell$ be the diffeomorphism defined by
\begin{eqnarray*}
{}&{}&H_5(X_1,X_2,\ldots,X_\ell)\\
{}&=&\left(X_1,X_2,\ldots,X_{\ell -1},X_\ell-
\sum_{k=1}^{\ell-1}\left(\frac{\sum_{j=\ell}^n\alpha_{j-\ell+1, k}^2}{\left(\sum_{k=1}^{\ell-1}\alpha_{jk}^2+1\right)}
+1\right)X_k^2
\right). 
\end{eqnarray*}
The composition of $H_5$ and $H_3\circ H_1\circ D_{(p_1,\ldots,p_\ell)}\circ H_2\circ H_4$ is as follows:
\begin{eqnarray*}
{}&{}&(H_5\circ H_3\circ H_1\circ D_{(p_1,\ldots,p_\ell)}\circ H_2\circ H_4)(x_1,x_2,\ldots,x_n)\\
{}&=&\left(x_1,x_2,\ldots,x_{\ell -1},
\sum_{j=\ell}^n\left(\sum_{k=1}^{\ell -1}\alpha_{jk}^2+1\right)
\left(x_j-\frac{1}{\left(\sum_{k=1}^{\ell -1}\alpha_{jk}^2+1\right)}\sum_{k=1}^{\ell-1}\alpha_{j-\ell+1, k}x_k\right)^2
\right).
\end{eqnarray*}
\par 
Let $H_6 : \R^n\rightarrow \R^n$ be the diffeomorphism defined by
\begin{eqnarray*}
{}&{}&H_6(x_1,x_2,\ldots,x_n)\\
{}&=&\biggl(x_1,\ldots,x_{\ell -1},
\frac{x_\ell}{\sqrt{r_\ell}}+\varphi_\ell(x_1,\ldots,x_{\ell-1}),\ldots,\frac{x_n}{\sqrt{r_n}}+\varphi _n(x_1,\ldots,x_{\ell-1})
\biggl), 
\end{eqnarray*}
where $r_j=\sum_{k=1}^{\ell -1}\alpha_{jk}^2+1$ and 
$\varphi_j(x_1,\ldots,x_{\ell-1})=\frac{1}{\left(\sum_{k=1}^{\ell -1}\alpha_{jk}^2+1\right)}\sum_{k=1}^{\ell-1}\alpha_{j-\ell+1, k}x_k$ $(\ell \le j\le n)$.    
It follows that
\begin{eqnarray*}
{}&{}&(H_5\circ H_3\circ H_1\circ D_{(p_1,\ldots,p_\ell)}\circ H_2\circ H_4\circ H_6)(x_1,x_2,\ldots,x_n)\\
{}&=&(x_1,x_2,\ldots,x_{\ell -1},x_\ell^2+\cdots +x_n^2).
\end{eqnarray*}
\hfill\qed 
\par 
\smallskip 
{\it Proof of $(II)$.}\qquad 
Since $n<\ell$ and $p_1, \ldots, p_{n+1}$ are in general position, 
there exists an $(\ell-1)\times (\ell-1)$ regular matrix $\widetilde{B}$ such that
\begin{eqnarray*}
A\widetilde{B}=
\begin{pmatrix}
1&&\mbox{\smash{\Large $0$}}&                                             0 & \cdots & 0   \\
&\ddots &&        \vdots &   & \vdots \\
\mbox{\smash{\Large $0$}}&&1&            0  &  \cdots  &  0\\
\end{pmatrix}, 
\end{eqnarray*}
where the matrix $A$ is the same as in the proof of $(I)$.     
Similarly as the proof of $(I)$ of Proposition \ref{P}, the composition $H_3\circ H_1\circ D_{(p_1,\ldots,p_\ell)}\circ H_2$ 
is as follows:
\begin{eqnarray*}
{}&{}&(H_3\circ H_1\circ D_{(p_1,\ldots,p_\ell)}\circ H_2)(x_1,x_2,\ldots,x_n)\\
{}&=&\begin{pmatrix}
x_1&x_2&\cdots &x_n
\end{pmatrix}
\left(
\begin{array}{ccc|c}
&&&  x_1\\ 
&&&x_2\\
&\mbox{\smash{\Large $A$}}&&\vdots \\
&&&\vdots\\
&&&x_n
\end{array}
\right)
\left(
\begin{array}{ccc|c}
&&&  0\\ 
&\mbox{\smash{\Large $\widetilde{B}$}}&&\vdots\\
&&&0\\ \hline
0&\cdots &0&1
\end{array}
\right)\\
{}&=&\begin{pmatrix}
x_1&x_2&\cdots &x_n
\end{pmatrix}
\begin{pmatrix}
1&&\mbox{\smash{\Large $0$}}&                                             0 & \cdots & 0  &x_1 \\
&\ddots &&        \vdots &   & \vdots & \vdots  \\
\mbox{\smash{\Large $0$}}&&1&            0  &  \cdots  &  0 & x_n\\
\end{pmatrix}\\
{}&=&\biggl(x_1,\ldots,x_n,0,\ldots,0,\sum_{j=1}^nx_j^2\biggl).
\end{eqnarray*}

Let $\widetilde{H}_4 : \R^\ell\rightarrow \R^\ell $ be the diffeomorphism defined by 
\begin{eqnarray*}
\widetilde{H}_4(X_1,X_2,\ldots,X_\ell)=\left(X_1,\ldots,X_n,\ldots,X_{\ell-1},X_\ell-\sum_{j=1}^nX_j^2\right).
\end{eqnarray*}
It follows that 
\begin{eqnarray*}
{}&{}&(\widetilde{H}_4\circ H_3\circ H_1\circ D_{(p_1,\ldots,p_\ell)}\circ H_2)(x_1,x_2,\ldots,x_n)\\
{}&=&(x_1,x_2,\ldots,x_n,0,\ldots,0).
\end{eqnarray*}
\hfill\qed 
\par 
\smallskip     
We have the following as a corollary of Proposition \ref{P} $(I)$. 
\begin{coro}\label{corollary 2.1}
In the case that $n=\ell$ $(\geq 2)$, let $p_1,\ldots,p_\ell\in \mathbb{R}^\ell$ $(p_{1\ell}=\cdots =p_{\ell\ell}=a)$ be in general position. Then, there exists a diffeomorphism $H:\mathbb{R}^\ell \to \mathbb{R}^\ell$ such that $H\circ D_{(p_1,\ldots , p_\ell)}$ 
is a definite fold mapping whose singular set is the hyperplane $\mathbb{R}^{\ell-1}\times\{ a\}$. 
\end{coro}
{\it Proof.}\qquad 
The composition $H_1\circ D_{(p_1,\ldots,p_\ell)}$ can be expressed as follows:
\begin{eqnarray*}
{}&{}&(H_1\circ D_{(p_1,\ldots,p_\ell)})(x_1,x_2,\ldots,x_\ell)\\
&=&\tr{
\begin{pmatrix}
x_1-p_{11}\\
\vdots \\ 
\vdots \\
x_\ell-a \\
\end{pmatrix}}
\begin{pmatrix}
p_{21}-p_{11}&\cdots &p_{\ell 1}-p_{11}&x_1-p_{11}\\
p_{22}-p_{12}&\cdots &p_{\ell 2}-p_{12}&x_2-p_{12}\\
\vdots & & \vdots &\vdots \\
p_{2,\ell-1}-p_{1,\ell-1}&\cdots &p_{\ell,\ell-1}-p_{1,\ell-1}&x_{\ell-1}-p_{1,\ell-1}\\
0 & \cdots & 0 &x_\ell-a
\end{pmatrix},
\end{eqnarray*}
where $\tr{\widetilde{A}}$ means the transposed matrix of $\widetilde{A}$.
The composition of $H_3$ and $H_1\circ D_{(p_1,\ldots,p_\ell)}$ is as follows:
\begin{eqnarray*}
{}&{}&(H_3\circ H_1\circ D_{(p_1,..., p_\ell)})(x_1,\,x_2,\,.\,.\,.\,,\,x_\ell)\\
{}&=&\biggl(x_1-p_{11},\ldots,x_{\ell-1}-p_{1,\ell-1},\sum_{j=1}^{\ell-1}(x_j-p_{1j})^2+(x_\ell-a)^2\biggl).
\end{eqnarray*}
Let $H_4' : \R^\ell\rightarrow \R^\ell$ be the diffeomorphism defined by
\[H_4'(X_1,X_2,\ldots,X_\ell)=\left(X_1+p_{11},\ldots,X_{\ell-1}+p_{1,\ell-1},X_\ell-\sum_{j=1}^{\ell-1}X_j^2+a\right).\]
Set $H=H_4'\circ H_3\circ H_1$, and we have the following:
\begin{eqnarray*}
{}&{}&(H\circ D_{(p_1,\ldots,p_\ell)})(x_1,x_2,\ldots,x_\ell)\\
{}&=&(x_1,x_2,\ldots,x_{\ell-1},(x_\ell-a)^2+a).
\end{eqnarray*}
\hfill\qed 
\par 
\smallskip
\begin{lemm}[\cite{mather, bruce}]\label{M}
Let $M$ be a closed manifold and let $i:M\rightarrow \R^\ell$ be an immersion with normal crossings.  
Then, there exists a subset $\Sigma\subset \mathcal{L}(\mathbb{R}^\ell,\mathbb{R})$ of Lebesgue measure $0$ such that for any $\pi \in \mathcal{L}(\mathbb{R}^\ell,\mathbb{R})-\Sigma$, the composition $\pi\circ i : M\rightarrow \mathbb{R}$ is a Morse function, where $\mathcal{L}(\mathbb{R}^\ell,\mathbb{R})$ stands for 
the space of linear functions $\R^\ell\rightarrow \R$.
\end{lemm}
Let $M$ be a closed manifold and let $i:M\rightarrow \R^\ell$ be an immersion with normal crossings.   
Then, by Lemma \ref{M} there exists a linear function 
$\pi : \mathbb{R}^\ell\rightarrow \mathbb{R}$: 
\[\pi (x)=\sum_{j=1}^{\ell}u_jx_j\hspace{5mm}(u=(u_1,\ldots,u_\ell),x=(x_1,\ldots,x_\ell))\]
such that $\pi\circ i : M\rightarrow \mathbb{R}$ is a Morse function.
By a change of basis of the linear space $\R^\ell$ if necessary, we may assume 
that  $u=(0,\ldots,0,1)$.
Note that $\pi (x)=x_\ell$, $\pi \circ i=i_\ell$ ($i=(i_1,\ldots,i_\ell)$).
Since $M$ is compact and $\pi\circ i$ is continuous, $\pi\circ i$ has the minimal value $k_0$.    
It is clear that $k_0$ is a critical value of $\pi\circ i$.
Since $\pi\circ i$ is a Morse function, the set $(\pi\circ i)^{-1}(k_0)$ consists of only one point $q$.
Let $(U,\varphi ,(t_1,\ldots,t_m))$ be a coordinate neighborhood of $q$. 
By shrinking $U$ sufficiently small if necessary, we may assume that $i|_{(U)}:U\rightarrow \mathbb{R}^\ell$ is an embedding. 

We show first that near $q$ a coordinate system can be chosen so that $i(M)$ can be expressed 
locally as a graph of a map $\mathbb{R}^m\to \mathbb{R}^{\ell-m}$.   
Let $(Ji)_q$ be the Jacobian matrix of the mapping $i$ at $q$. Note that 
\[\frac{\partial (i_\ell\circ \varphi ^{-1})}{\partial t_j}(\varphi (q))=0\]
for any $j$ ($1\leq j\leq m$).
Thus, we have
\[(0,\ldots,0,1)(Ji)_q=(0,\ldots,0).\]
By a change of basis of the linear space $\R^\ell$ if necessary, 
we may assume 
that the matrix 
\[
(E_m\mid O)(Ji)_q
\]
is invertible, 
where $E_m$ is the $m\times m$ unit matrix and $O$ is the zero $m\times (\ell-m)$ matrix.

Consider the mapping $\psi : \varphi (U)\rightarrow \R^m\times \{0\}\subset \R^\ell$ defined by
\[
\psi =(i_1\circ \varphi ^{-1},\ldots,i_m\circ \varphi ^{-1}).\]
It follows that $(J\psi )_{\varphi (q)}$ is invertible.
By the inverse function theorem, there exist an open neighborhood $V$ of $\varphi (q)$ and 
an open neighborhood $\widetilde {V}$ of $\psi$ ($\varphi (q)$) such that $\psi : V\rightarrow \widetilde {V}$ is 
a diffeomorphism. We may suppose that $V$ and $\widetilde {V}$ are connected and $V\subset \varphi (U)$.
Then, $i(\varphi ^{-1}(V))$ can be expressed as follows:
\[i(\varphi ^{-1}(V))=\{(x,\xi _{m+1}(x),\ldots,\xi _\ell(x))\in \mathbb{R}^\ell
\mid x\in \widetilde {V}\},\]
where $x=(x_1,\ldots,x_m)$ and $\xi _j=i_j\circ \varphi ^{-1}\circ \psi ^{-1}$ for any $j$ $(m+1\leq j\leq \ell)$.
\par 
Set 
\[
k_1=\displaystyle\sup_{q'\in \varphi ^{-1}(V)}i_\ell(q').
\]
Then, $k_1>k_0$, and $i_\ell^{-1}(r)$ is a non-empty set for any $r\in [k_0,k_1)$.
Since $M-\varphi ^{-1}(V)$ is compact and $i_\ell$ is continuous, 
there exists a minimal value 
$k_2$ of $i_\ell$.    
This value satisfies $k_2>k_0$.

Set 
\[
k_3=\min\Bigl\{k_0+1,k_0+\displaystyle\frac{k_2-k_0}{3},k_1\Bigr\}.
\]
\begin{lemm}\label{lemma 2.1}
There exists a real number $a\in (k_0,k_3)$ such that 
there exist $(m+1)$-points in general position in 
$i(M)\cap (\mathbb{R}^{\ell-1}\times \{a\})$.   
\end{lemm} 
{\it Proof.}\qquad 
Suppose that $m\geq 2$. By Sard's theorem, there exists a real number $a\in (k_0,k_3)$ such that $a$ is a regular value of the mapping $i_\ell : M\rightarrow \R$. 
Since a real number $a$ is also an element of $(k_0,k_1)$, we get $i_\ell^{-1}(a)\not=\emptyset$.
Hence, $i_\ell^{-1}(a)$ is an $(m-1)$-dimensional closed manifold.
Note that $i(i_\ell^{-1}(a))=i(M)\cap (\R^{\ell-1}\times \{a\})$. Let $n_0$ be the maximal value of the number 
of points which are in general position in $i(i_\ell^{-1}(a))$. 
It is clear that $n_0\geq 2$. 
Let $p_1,\ldots,p_{n_0}$ be $n_0$-points in general position in $i(i_\ell^{-1}(a))$, and set 
\[
W=\Bigl\{\displaystyle\sum_{j=2}^{n_0}\alpha _j\overrightarrow{p_1p_j}\mid \alpha _j\in \R\Bigr\}.
\]
Then, $i(i_\ell^{-1}(a))\subset W$ and ${\rm dim}\;W=n_0-1$. 
Since $i\mid _{i_\ell^{-1}(a)} : i_\ell^{-1}(a)\rightarrow W$ is an immersion and 
${i_\ell^{-1}(a)}$ is closed, it follows that $m-1<n_0-1$.

Suppose that $m=1$.
Since $\widetilde {V}$ is a bounded connected open set of $\R$, we can set $\widetilde {V}=(b_1,b_2)$.
Then, the set $i(\varphi ^{-1}(V))$ can be expressed as follows:
\[i(\varphi ^{-1}(V))=\{(x,\xi _2(x),\ldots,\xi _\ell(x))\in \R^\ell \mid x\in (b_1,b_2)\}.\]
Note that $b_1<i_1(q)<b_2$.    Define the following two closed intervals:
\begin{eqnarray*}	
B_1&=&\Bigl[b_1+\frac{i_1(q)-b_1}{2},\,i_1(q)\Bigr],\\
B_2&=&\Bigl[i_1(q),\,i_1(q)+\frac{b_2-i_1(q)}{2}\Bigr].
\end{eqnarray*}
Then, 
\begin{eqnarray*}	
\xi _\ell\Bigl(b_1+\frac{i_1(q)-b_1}{2}\Bigr)&>&\xi _\ell(i_1(q))=k_0,\\
\xi _\ell\Bigl(i_1(q)+\frac{b_2-i_1(q)}{2}\Bigr)&>&\xi _\ell(i_1(q))=k_0.
\end{eqnarray*}
We set  
\begin{eqnarray*}	
a=\min\Bigl\{\xi _\ell\Bigl(b_1+\frac{i_1(q)-b_1}{2}\Bigr),\,\xi _\ell\Bigl(i_1(q)+\frac{b_2-i_1(q)}{2}\Bigr),\,\frac{k_0+k_3}{2}\Bigr\}.
\end{eqnarray*}
It follows that $a\in (k_0,k_3)$. 
Since $B_1$ and $B_2$ are compact and $\xi _\ell$ is continuous, by the intermediate value theorem, there exists 
a real number $c_1\in B_1$ such that $\xi _\ell(c_1)=a$ and there exists a real number $c_2\in B_2$ such that $\xi _\ell(c_2)=a$. Since $c_1\not=c_2$, there exist at least two points 
$(c_1,\xi _2(c_1),\ldots,\xi _\ell(c_1))$, $(c_2,\xi _2(c_2),\ldots,\xi _\ell(c_2))$ 
in $i(M)\cap (\R^{\ell-1}\times \{a\})$.
\hfill\qed    
\par 
\smallskip 
\subsection{Preliminaries for the proof of Proposition \ref{d}}   
The proof of Proposition \ref{d} requires 
the following four lemmas. 
\begin{lemm}[\cite{Berger}, p.354]\label{AH}
Let $A\subset \R^2$ be a closed convex set, and let $x$ be a point in the boundary $\partial A$ of $A$.  
Then, there exist the following sets
\begin{eqnarray*}	
L=\{(x_1,x_2)\in \R^2\mid \alpha x_1+\beta x_2+\gamma =0\},\\
L_{\leq }=\{(x_1,x_2)\in \R^2\mid \alpha x_1+\beta x_2+\gamma \leq 0\}
\end{eqnarray*}
such that $x\in L$, $A\subset L_{\leq }$ $(\alpha \not=0 \mbox{ or } \beta \not=0)$.
\end{lemm}
\begin{lemm}[\cite{Berger}, p.339]\label{AK}
Let $A$ be a subset of $\R^2$. Then, for any point $x=(x_1,x_2)$ in the convex hull $\mathrm{conv}(A)$ of $A$,  there exist three points $q_1=(q_{11},q_{12})$, $q_2=(q_{21},q_{22})$, $q_3=(q_{31},q_{32})$ $\in A$ 
and real numbers $t_1$,$t_2$,$t_3$ such that 
\begin{eqnarray*}	
x_1&=&t_1q_{11}+t_2q_{21}+t_3q_{31}\\
x_2&=&t_1q_{12}+t_2q_{22}+t_3q_{32}\\
&{}&(t_1\geq 0,\,t_2\geq 0,\,t_3\geq 0,\,t_1+t_2+t_3=1).   
\end{eqnarray*}
\end{lemm}
\begin{lemm}[\cite{Berger}, p.343]\label{AS}
Let $A$ be a subset of $\R^2$ homeomorphic to $S^1$. Then, $\partial \mathrm{conv}(A)$ is homeomorphic to $S ^1$.
\end{lemm}
Now, we define the following mapping which is important in the proof of Proposition \ref{d}.
\begin{defi}\label{definition 2.1}
{\rm Let $\pi _\theta : i(S^1)\rightarrow \R$ be the mapping defined by 
\[\pi _\theta (x_1,x_2)=(\cos\theta )x_1+(\sin\theta )x_2.\]}
\end{defi}
Since $i(S^1)$ is compact and $\pi _\theta $ is continuous, 
$\pi _\theta $ has the maximal value $k_\theta$.   
Then, $(I)$ or $(II)$ holds.   
\medskip \\
$(I)$ For any $\theta \in [0,2\pi )$, the set $\pi _\theta ^{-1}(k_\theta )$ has only one point.\\
$(II)$ There exists $\theta \in [0,2\pi )$ such that the set $\pi _\theta ^{-1}(k_\theta )$ has at least two points.   
\begin{lemm}\label{I}
In the circumstance of $(I)$, $i(S^1)=\partial \mathrm{conv}(i(S^1))$.
\end{lemm}
{\it Proof.}\qquad  
First, we prove that $\partial \mathrm{conv}(i(S^1))\subset i(S^1)$. 
Suppose that there exists $a\in$ $\partial \mathrm{conv}(i(S^1))$ such that $a\not\in i(S^1)$ $(a=(a_1,a_2))$. 
Then, by Lemma \ref{AH}, there exist the following sets: 
\begin{eqnarray*}	
L=\{(x_1,x_2)\in \R^2\mid \alpha x_1+\beta x_2+\gamma =0\},\\
L_{\leq }=\{(x_1,x_2)\in \R^2\mid \alpha x_1+\beta x_2+\gamma \leq 0\}
\end{eqnarray*}
such that $a\in L$, $\mathrm{conv}(i(S^1))\subset L_{\leq }$. 
By changing the origin and basis of the linear space $\R^2$ if necessary, we may assume that 
\begin{eqnarray*}	
a&=&(0,0),\\
L&=&\R\times \{0\},\\
L_{\leq }&=&\R\times \{x\in \R\mid x\leq 0\}.
\end{eqnarray*}

Note that $(a_1,a_2)=(0,0)$, and by Lemma \ref{AK} there exist three points $q_1=(q_{11},q_{12})$, 
$q_2=(q_{21},q_{22})$, $q_3=(q_{31},q_{32})$ $\in i(S^1)$ and 
real numbers $t_1$,$t_2$,$t_3$ such that 
\begin{eqnarray}	
&&t_1q_{11}+t_2q_{21}+t_3q_{31}=0\\
&&t_1q_{12}+t_2q_{22}+t_3q_{32}=0\\
&&(t_1\geq 0,\,t_2\geq 0,\,t_3\geq 0,\,t_1+t_2+t_3=1).\nonumber
\end{eqnarray}
Note that $q_{12}\leq 0$, $q_{22}\leq 0$, $q_{32}\leq 0$, and 
there exists $j$ $(1\leq j\leq 3)$ such that $q_{j2}=0$. 
Then, we may assume that $q_{12}=0$. 

Thus, the maximal value of the mapping $\pi _{\pi /2}$ is $0$. 
Since the set $\pi _{\pi /2}^{-1}(0)$ has only one point, it follows that $q_{22}<0$, $q_{32}<0$. 
Hence, by $(2)$, it is necessary that $t_2=t_3=0$. 
Since we have $t_1=1$, we obtain $a=q_1$ by $(1)$ and $(2)$. This is inconsistent with $a\not\in i(S^1)$.

Next, we prove that $i(S^1)\subset \partial \mathrm{conv}(i(S^1))$. Suppose that there exists $a\in i(S^1)$ such that 
$a\not\in \partial \mathrm{conv}(i(S^1))$. Since $\partial \mathrm{conv}(i(S^1))\subset i(S^1)$, we obtain 
$\partial \mathrm{conv}(i(S^1))\subset i(S^1)-\{a\}$. 
This means that $\partial \mathrm{conv}(i(S^1))$ is not homeomorphic to $S^1$.
This is inconsistent with Lemma \ref{AS}. \qed 
\section{Proofs of Theorems \ref{D1} and \ref{D2}}\label{SD}
Except the last parts, the proof of Theorem \ref{D1} is completely 
the same as the proof of Theorem \ref{D2} as follows.   
\par 
By Lemma \ref{M} there exists a linear function 
$\pi : \mathbb{R}^\ell\rightarrow \mathbb{R}$: 
\[\pi (x)=\sum_{j=1}^{\ell}u_jx_j\hspace{5mm}(u=(u_1,\ldots,u_\ell),x=(x_1,\ldots,x_\ell))\]
such that $\pi\circ i : M\rightarrow \mathbb{R}$ is a Morse function.
Let $k_0$ be the minimal value of $\pi\circ i$ and let $q\in M$ be the unique point satisfying 
$(\pi\circ i)(q)=k_0$ as in Section 2.       
Then, by Lemma \ref{lemma 2.1}, 
there exists a real number $a\in (k_0,k_3)$ such that there exist $(m+1)$-points $p_1,\ldots,p_{m+1}$ 
in general position in 
$i(M)\cap (\R^{\ell-1}\times \{a\})$. 
Then, we take $\ell-(m+1)$-points \;$p_{m+2},\ldots,p_\ell\in \mathbb{R}^{\ell-1}\times \{a\}$\;  
such that $\ell$-points $p_1,\ldots,p_\ell$ are in general position.

Note that $p_{1\ell}=\cdots=p_{\ell\ell}=a$.   
Thus, by Corollary \ref{corollary 2.1}, the mapping $H\circ D_{(p_1,\ldots,p_\ell)}$ is a definite fold mapping whose singular set is the hyperplane 
$\R^{\ell-1}\times \{a\}$. Hence, any singular point of $(H\circ D_{(p_1,\ldots,p_\ell)})\mid _{i(M)}$ must be 
contained in $i(M)\cap (\R^{\ell-1}\times \{a\})$. Moreover, since the intersection $i(M)\cap (\R^{\ell-1}\times \{a\})$ 
is contained in $i(\varphi ^{-1}(V))$, 
the following equalities hold for any point in $i(M)\cap (\R^{\ell-1}\times \{a\})$:
\begin{eqnarray*}	
{}&{}&(H\circ D_{(p_1,\ldots,p_\ell)})\mid _{i(M)}(x,x_{m+1},\ldots,x_\ell)\\
{}&=&(H\circ D_{(p_1,\ldots,p_\ell)})\mid _{i(M)}(x,\xi _{m+1}(x),\ldots,\xi _\ell(x))\\
{}&=&(x,\xi _{m+1}(x),\ldots,\xi _{\ell-1}(x),(\xi _\ell(x)-a)^2+a),
\end{eqnarray*}
where $x=(x_1,\ldots,x_m)$. 
Note that the rank of the Jacobian matrix of $(H\circ D_{(p_1,\ldots,p_\ell)})\mid _{i(M)}$ at any point of $i(M)\cap (\R^{\ell-1}\times \{a\})$ is always $m$, which implies that $H\circ D_{(p_1,\ldots,p_\ell)}\circ i$ is non-singular.   
\par  
Next, we show that $(H\circ D_{(p_1,\ldots,p_\ell)})\mid _{i(M)}$ is injective.     
For a real number $r$, set 
\begin{eqnarray*}	
A_{<r}&=&i(M)\cap (\mathbb{R}^{\ell-1}\times \{x\in \mathbb{R}\mid x<r\}),\\
A_{\geq r}&=&i(M)\cap (\mathbb{R}^{\ell-1}\times \{x\in \mathbb{R}\mid x\geq r\}).
\end{eqnarray*}
Since $H\circ D_{(p_1,\ldots,p_\ell)}$ is a definite fold mapping whose singular set is the hyperplane $\R^{\ell-1}\times \{a\}$ by Corollary \ref{corollary 2.1}, the following holds.    
\[ 
(H\circ D_{(p_1,\ldots,p_\ell)})\mid _{i(M)} \mbox{ is injective }
\Leftrightarrow  (H\circ D_{(p_1,\ldots,p_\ell)})(A_{<a})\cap A_{\geq a}=\emptyset.
\]

In order to prove $(H\circ D_{(p_1,\ldots,p_\ell)})(A_{<a})\cap A_{\geq a}=\emptyset$, it is sufficient to show the following $(3)$ and $(4)$:
\begin{eqnarray}
(H\circ D_{(p_1,\ldots,p_\ell)})(A_{<a})\cap A_{\geq k_2}&=&\emptyset,  \\
(H\circ D_{(p_1,\ldots,p_\ell)})(A_{<a})\cap A_{<k_2}&=&\emptyset.   
\end{eqnarray}

First, we prove $(3)$. Suppose that $(H\circ D_{(p_1,\ldots,p_\ell)})(A_{<a})\cap A_{\geq k_2}\not=\emptyset $.
Then, there exists a real number $\widetilde {a}\in [k_0,a)$ such that $(\widetilde {a}-a)^2+a\geq k_2$.
Then, we have $(k_0-a)^2+a\geq (\widetilde {a}-a)^2+a$.  It follows that $(k_0-a)^2+a\geq k_2$.
On the other hand, since $k_0<a<k_0+1$, we obtain $0<a-k_0<1$. 
Thus, we have $(k_0-a)^2+a<2(a-k_0)+a$. Since $a<k_0+(k_2-k_0)/3$, it follows that $2(a-k_0)+a<k_2$. 
Therefore, we obtain $(k_0-a)^2+a<k_2$. However, this is inconsistent with $(k_0-a)^2+a\geq k_2$.

For $(4)$,  
suppose that $(H\circ D_{(p_1,\ldots,p_\ell)})(A_{<a})\cap A_{<k_2}\not=\emptyset$.
Note that $A_{<k_2}\subset i(\varphi ^{-1}(V))$, and 
there exist two points \[(x,\xi _{m+1}(x),\ldots,\xi _\ell(x))\in A_{<a},\ (x',\xi _{m+1}(x'),\ldots,\xi _\ell(x'))\in A_{<k_2}\]
such that 
\begin{eqnarray*}	
(H\circ D_{(p_1,\ldots,p_\ell)})(x,\xi _{m+1}(x),\ldots,\xi _\ell(x))=(x',\xi _{m+1}(x'),\ldots,\xi _\ell(x')),
\end{eqnarray*}
where $x=(x_1,\ldots,x_m)$ and $x'=(x_1',\ldots,x_m')$. 
Since $x=x'$ and $\xi _\ell(x')=(\xi_\ell(x)-a)^2+a$, we get $\xi _\ell(x)\geq a$. However, this contradicts to 
$\xi _\ell(x)<a$.     
Therefore, $(H\circ D_{(p_1,\ldots,p_\ell)})\mid _{i(M)}$ must be injective. 
\subsection{Proof of Theorem \ref{D1}}\label{D1P}
We have already proved that $H\circ D_{(p_1,\ldots,p_\ell)}\circ i$ is an immersion and 
$(H\circ D_{(p_1,\ldots,p_\ell)})\mid _{i(M)}$ is injective.        
Since $i$ is an embedding and $H$ is a diffeomorphism,  
the mapping $D_{(p_1,\ldots,p_\ell)}\circ i$ must be an embedding.  
\hfill\qed    
\par 
\subsection{Proof of Theorem \ref{D2}}
It is sufficient to show that  $H\circ D_{(p_1,\ldots,p_\ell)}\circ i$ is with normal crossings. 
Since $(H\circ D_{(p_1,\ldots,p_\ell)})\mid _{i(M)}$ is injective, 
for any $y\in (H\circ D_{(p_1,\ldots,p_\ell)}\circ i)(M)$, the set $((H\circ D_{(p_1,\ldots,p_\ell)})\mid _{i(M)})^{-1}(y)$ 
has only one point. 
Let $x$ be that point. 

Suppose that $x\in i(M)\cap (\R^{\ell-1}\times \{a\})$. Since $i^{-1}(x)\subset \varphi ^{-1}(V)$ and 
$i\mid _{\varphi ^{-1}(V)} :  \varphi ^{-1}(V) \rightarrow \R^\ell$ is an embedding, the set $i^{-1}(x)$ has 
only one point.     
It follows that $H\circ D_{(p_1,\ldots,p_\ell)}\circ i$ is a mapping with normal crossings at the 
point.

Suppose that $x\in i(M)-(\R^{\ell-1}\times \{a\})$. Since $i$ is an immersion with normal crossings, 
the set $(H\circ D_{(p_1,\ldots,p_\ell)}\circ i)^{-1}(y)$ can be expressed as follows: 
\[(H\circ D_{(p_1,\ldots,p_\ell)}\circ i)^{-1}(y)=i^{-1}(x)=\{q_1,\ldots,q_n\}.\]
Since $d(H\circ D_{(p_1,\ldots,p_\ell)})_x : T_x\R^\ell\rightarrow T_y\R^\ell$ is a bijective linear mapping, for any subset 
$\{\lambda _1,\lambda _2,\ldots,\lambda _s\}\subset \{1,2,\ldots,n\}$ $(s\leq n)$, the following equalities hold.   
\begin{eqnarray*}	
{}&{}&{\rm codim}\,\Bigl(\bigcap _{j=1}^sd(H\circ D_{(p_1,\ldots,p_\ell)}\circ i)_{q_{\lambda _{j}}}(T_{q_{\lambda _{j}}}M)\Bigr)\\
{}&=&\ell-{\rm dim}\,\bigcap _{j=1}^sd(H\circ D_{(p_1,\ldots,p_\ell)}\circ i)_{q_{\lambda _{j}}}(T_{q_{\lambda _{j}}}M)\\
{}&=&\ell-{\rm dim}\,d(H\circ D_{(p_1,\ldots,p_\ell)})_x\Bigl(\bigcap _{j=1}^sdi_{q_{\lambda _{j}}}(T_{q_{\lambda _{j}}}M)\Bigr)\\
&=&\ell-{\rm dim}\,\Bigl(\bigcap _{j=1}^sdi_{q_{\lambda _{j}}}(T_{q_{\lambda _{j}}}M)\Bigr)\\
&=&\sum_{j=1}^s(\ell-{\rm dim}\,di_{q_{\lambda _{j}}}(T_{q_{\lambda _{j}}}M))\\
&=&\sum_{j=1}^s(\ell-{\rm dim}\,d(H\circ D_{(p_1,\ldots,p_\ell)}\circ i)_{q_{\lambda _{j}}}(T_{q_{\lambda _{j}}}M))\\
&=&\sum_{j=1}^s{\rm codim}\,d(H\circ D_{(p_1,\ldots,p_\ell)}\circ i)_{q_{\lambda _{j}}}(T_{q_{\lambda _{j}}}M)
\end{eqnarray*}
Hence, $H\circ D_{(p_1,\ldots,p_\ell)}\circ i$ is an immersion with normal crossings. 
\hfill\qed 
\section{Proof of Proposition \ref{d}}\label{dP}
Since $i(S^1)$ is compact, $\R^2$ is a Hausdorff space, $d_{(p_1,p_2)}\mid _{i(S^1)}$ is continuous 
and $i$ is homeomorphic to the image, 
in order to prove that $d_{(p_1,p_2)}\circ i$ is homeomorphic to the image, it is sufficient to show that 
$d_{(p_1,p_2)}\mid _{i(S^1)}$ is injective.    
Recall that $k_\theta$ is the maximal value of the map $\pi_\theta$  defined in Definition \ref{definition 2.1} and  
one of the following $(I)$ and $(II)$ holds for $k_\theta$.   
\par 
\smallskip 
\noindent 
$(I)$ For any $\theta \in [0,2\pi )$, the set $\pi _\theta ^{-1}(k_\theta )$ has only one point.\\
$(II)$ There exists $\theta \in [0,2\pi )$ such that the set $\pi _\theta ^{-1}(k_\theta )$ has at least two points.   
\par 
We consider first the case $(I)$.     
Let $d : \R^2\times \R^2\rightarrow \R$ be the two-dimensional Euclidean distance defined in Section \ref{IN}.
Since the set $i(S^1)\times i(S^1)$ is compact and $d$ is continuous, $d\mid _{i(S^1)\times i(S^1)}$ 
has a maximal value $k$.   
Let $(a,a')$ $\in i(S^1)\times i(S^1)$ be an element of the set $(d\mid _{i(S^1)\times i(S^1)})^{-1}(k)$.
By changing the origin and basis of the linear space $\R^2$ if necessary, 
we may assume  
that $a=(0,0)$, $a'=(1,0)$. Then, $k=1$.

Set $\theta_0 =\pi /2$, and define $\pi _{\theta _0} : i(S^1)\rightarrow \R$ ; thus 
\begin{eqnarray*}	
\pi _{\theta _0}(x_1,x_2)=x_2.
\end{eqnarray*}
Since $i(S^1)$ is compact and $\pi _{\theta _0}$ is continuous, $\pi _{\theta _0}$ has the maximal value 
$k_{\theta _0}$ and the minimal value $\widetilde k_{\theta _0}$.    
Since there exist two points $a=(0,0)$, $a'=(1,0)$ in $i(S^1)$, it is clear that 
$\widetilde k_{\theta _0}\leq 0\leq k_{\theta _0}$. 
Then, the set $\pi _{\theta _0}^{-1}(k_{\theta _0})$ has only one point. Since $-\widetilde k_{\theta _0}$ is 
also the maximal value of $\pi_{3\pi/2}$ and $\pi _{\theta _0}^{-1}(\widetilde k_{\theta _0})=\pi_{3\pi/2}^{-1}(-\widetilde k_{\theta _0})$, the set $\pi _{\theta _0}^{-1}(\widetilde k_{\theta _0})$ 
also has only one point.    
It follows that 
$\widetilde k_{\theta _0}<0<k_{\theta _0}$. 
\par 
Again by changing of the basis of the linear space $\R^2$ if necessary, we may assume  
that $\mid \widetilde k_{\theta _0}\mid \leq \mid k_{\theta _0}\mid $.

We shall prove that the set $\pi _{\theta _0}^{-1}(2k_{\theta _0}/3)$ has at least two points. 
Note that $\pi _{\theta _0}^{-1}(2k_{\theta _0}/3)=i(S^1)\cap (\R\times \{2k_{\theta _0}/3\})$. Suppose that $\pi _{\theta _0}^{-1}(2k_{\theta _0}/3)=\emptyset $.    
Since $i(S^1)\cap (\R\times \{k_{\theta _0}\})\not=\emptyset$ and $i(S^1)\cap (\R\times \{0\})\not=\emptyset $, 
this is inconsistent with the fact that $i(S^1)$ is connected. Suppose that the set $\pi _{\theta _0}^{-1}(2k_{\theta _0}/3)$ has only one point, and 
put $\pi _{\theta _0}^{-1}(2k_{\theta _0}/3)=\{q\}$.    
Since $i(S^1)-\{q\}$ is not connected in this case, an contradiction yields.
Hence, the set $\pi _{\theta _0}^{-1}(2k_{\theta _0}/3)$ has at least two points. 
Let $p_1$,$p_2$ be the two points such that 
\begin{eqnarray*}	
p_1=\Bigl(p_{11},\frac{2}{3}k_{\theta _0}\Bigr),\ p_2=\Bigl(p_{21},\frac{2}{3}k_{\theta _0}\Bigr)\in \pi _{\theta _0}^{-1}(2k_{\theta _0}/3)\ (p_{11}<p_{21}).
\end{eqnarray*}

Suppose that there exist two points $b$, $b'\in i(S^1)$ $(b\not=b')$ such that 
\[d_{(p_1,p_2)}\mid _{i(S^1)}(b)=d_{(p_1,p_2)}\mid _{i(S^1)}(b'),\leqno{(5)}\]
where $b=(b_1,b_2)$ and $b'=(b'_1,b'_2)$. Note that $b\not=b'$, and by $(5)$ we obtain 
\begin{eqnarray*}	
b'_1&=&b_1,\\
b'_2&=&-b_2+\frac{4}{3}k_{\theta _0}.
\end{eqnarray*}
Since $b\not=b'$ and $b_1=b'_1$, we may assume that $b_2<b'_2$. Now, we show that $0<b_1<1$.
If $b_1>1$ (resp., $b_1<0$), we have $d\mid _{i(S^1)\times i(S^1)}(a,b)>1$ (resp., $d\mid _{i(S^1)\times i(S^1)}(a',b)>1$).   
These are contradictory with the fact that the maximal value of $d\mid _{i(S^1)\times i(S^1)}$ 
is $1$. 
If $b_1=0$, there exist two points $b$, $b'$ in $\pi _\pi ^{-1}(0)$ and if $b_1=1$, there exist two points 
$b$, $b'$ in $\pi _0^{-1}(1)$. Since the maximal value of $\pi _\pi$ is $0$ and  the maximal value of $\pi _0$ is $1$, these are inconsistent with $(I)$. Hence, we have $0<b_1<1$. 
Since $b_2<4k_{\theta _0}/3$ (resp., $b'_2<4k_{\theta _0}/3$), by 
$b'_2=-b_2+\frac{4}{3}k_{\theta _0}$,  
we have $b'_2>0$ (resp., $b_2>0$). 

Set 
\begin{eqnarray*}	
\bigtriangleup aa'b'=\{t_1a+t_2a'+t_3b'\in \R^2\mid t_1,t_2,t_3\geq 0,\ t_1+t_2+t_3=1\}.
\end{eqnarray*}

\begin{figure}
\begin{center}
\includegraphics[width=6cm,clip]{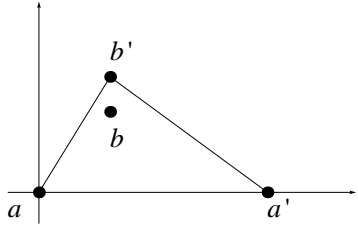}
\end{center}
\caption{$\bigtriangleup aa'b'$ and the point $b$}
\end{figure}
Since the set $i(S^1)$ has three points $a$,$a'$,$b'$ and by Lemma \ref{I} $i(S^1)=\partial \mathrm{conv}(i(S^1))$, 
it is clearly seen that $\bigtriangleup aa'b'\subset \mathrm{conv}(i(S^1))$. 
Then, $(\bigtriangleup aa'b')^\circ \subset (\mathrm{conv}(i(S^1)))^\circ $, where $A^\circ \subset \R^2$ 
is the set which consists of the interior points of $A$. 
It follows that $b\in (\bigtriangleup aa'b')^\circ $ (see Figure 2).    
Hence, we have $b\in (\mathrm{conv}(i(S^1)))^\circ $. 
This is inconsistent with $b\in \partial \mathrm{conv}(i(S^1))$.
\par 
\smallskip 
We consider now the case $(II)$. 
There exists $\theta _1\in [0,2\pi )$ such that the set $\pi _{\theta _1}^{-1}(k_{\theta _1})$ has at least 
two points, where $k_{\theta _1}$ is the maximal value of $\pi_{\theta _1}$. 
Let $p_1$,$p_2$ be two points in $\pi _{\theta _1}^{-1}(k_{\theta _1})$.   
We may assume 
that $p_1=(0,0)$, $p_2=(1,0)$, and $i(S^1)\subset \R\times \{x\in \R\mid x\leq 0\}$. 
Then, note that $\theta _1=\pi /2$, $k_{\theta_1}=0$. 
Now, suppose that there exist $c=(c_1,c_2)$, $c'=(c'_1,c'_2)\in i(S^1)$ $(c\not=c')$ such that 
\[d_{(p_1,p_2)}\mid _{i(S^1)}(c)=d_{(p_1,p_2)}\mid _{i(S^1)}(c').\leqno{(6)}\]

Note that $c\not=c'$, and by $(6)$, we obtain 
\begin{eqnarray*}	
c_1&=&c'_1,\\
c_2&=&-c'_2.
\end{eqnarray*}
Since $c_2\leq 0$ and $c'_2\leq 0$, it follows that $c_2=c'_2=0$.    This contradicts to the assumption  $c\not=c'$. 
\hfill\qed
\section*{Acknowledgement}
The authors are most grateful to the anonymous reviewer for his/her careful reading 
of the first manuscript of this paper and invaluable suggestions.    
They would like, also, to thank O. Saeki and M. Nishioka for their valuable comments.      
T. N. was partially supported 
by JSPS and CAPES under the Japan--Brazil research cooperative 
program.     




\bibliographystyle{model1-num-names}



\end{document}